\documentclass{amsart}
\usepackage{amsmath, amssymb, a4}

\numberwithin{equation}{section}
\theoremstyle{plain}
\newtheorem{thm}{Theorem}[section]
\newtheorem{lem}[thm]{Lemma}

\newtheorem{cor}[thm]{Corollary}

\theoremstyle{definition}

\newcommand{\q}{\quad}

\newcommand{\mrm}{\mathrm}

\newcommand{\sltwo}{\mrm{SL}(2, \mathbb{Z})}

\newcommand{\thalf}{\tfrac{1}{2}}
\newcommand{\summ}{\mathop{{\sum}^{\star}}}

\newcommand{\intt}{\int_{-\infty}^{\infty}}

\begin{document}

\title[Simultaneous nonvanishing of $L$-functions]{Simultaneous nonvanishing of Dirichlet $L$-functions and twists of Hecke-Maass L-functions}

\author{Soumya Das} 
\address{Department of Mathematics\\
Indian Institute of Science\\
Bangalore, India.}
\email{somu@math.iisc.ernet.in}

\author{Rizwanur Khan}
\address{
Science Program\\ Texas A\&M University at Qatar\\ Doha, Qatar}
\email{rizwanur.khan@qatar.tamu.edu}


\begin{abstract}
We prove that given a Hecke-Maass form $f$ for $\sltwo$ and a sufficiently large prime $q$, there exists a primitive Dirichlet character $\chi$ of conductor $q$ such that the $L$-values $L(\thalf, f \otimes \chi )$ and $L(\thalf, \chi)$ do not vanish. We expect the same method to work for any large integer $q$.
\end{abstract}
\maketitle

\section{Introduction}
Understanding when an $L$-function's central value is non-zero is a problem of great significance and long history in analytic number theory. Understanding when two or more $L$-functions are simultaneously non-zero is also a popular topic and can have important applications. For example, the problems considered in \cite{is} and \cite{michel} have a bearing on the Landau-Siegel zero problem and the Birch-Swinnerton-Dyer conjecture respectively. Recently there has been some interest in the simultaneous non-vanishing of $L$-functions in the family of primitive Dirichlet characters. Blomer and Mili\'cevi\'c \cite{blomil} proved that given two fixed Hecke-Maass forms $f_1$ and $f_2$ which satisfy the Ramanujan-Petersson conjecture, and a sufficiently large integer $r$ subject to a technical condition that does not permit integers such as primes and the product of two primes of almost equal size, there exists a primitive Dirichlet character $\chi$ of modulus $r$ such that $L(1/2, f_1\otimes \chi)$ and $L(1/2, f_2\otimes \chi)$ are not zero. This improved the work of Hoffstein and Lee \cite{hoflee}, who studied the same kind of non-vanishing problem, but were not able to specify the modulus of the character. We study a simpler problem: that of the simultaneous non-vanishing $L(1/2, f \otimes \chi)$ and $L(1/2, \chi)$ for a fixed Hecke-Maass form $f$ for $\sltwo$.

Our paper is motivated by a recent result of Liu \cite{liu}, who showed that for $R$ large enough, there exists a Dirichlet character modulo $r$ with $R<r<4R$ such that $L(1/2, f \otimes \chi)$ and $L(1/2, \chi)$ are non-vanishing. Actually Liu considers more general automorphic forms, but the point is that the modulus of the desired Dirichlet character is not specified by his work. Liu's result follows by showing that
\begin{align}
\sum_{r\in \mathfrak{D}} \summ_{\substack{\chi \bmod r\\ \chi(-1)=1}} L(1/2, f \otimes \chi) L(1/2, \chi) =  \frac{1}{2} \sum_{r\in \mathfrak{D}}  r + O(R^{\frac{15}{8}+\epsilon})
\end{align}
for some subset $\mathfrak{D}$ of the integers between $R$ and $4R$ of size $|\mathfrak{D}| \gg R/(\log R)^2$, where $\star$ indicates that the sum is restricted to the primitive characters. Liu notes that since the right hand side is non-zero, at least one of the summands on the left hand side must be non-zero. The main difficulty in Liu's analysis seems to be in dealing with Gauss sums, which arise from approximate functional equations. 

We will prove the following result.

\begin{thm} \label{mainthm} \label{mainthm}
Let $f$ be a Hecke-Maass form for $\sltwo$. For prime values of $q$, we have that
\begin{align}
\label{thmline} \summ_{\substack{\chi \bmod q \\ \chi(-1)=1}} L(\thalf, f \otimes \chi ) \overline{L(\thalf, \chi)} = \frac{q-2} {2} L(1,f) + O(q^{\frac{7}{8} + \theta + \epsilon}),
\end{align}
where $\epsilon>0$ is arbitrarily small. The implied constant depends on $f$ and $\epsilon$. In the exponent, $\theta$ represents the best bound towards the Ramanujan-Petersson conjecture for $f$, which can currently be taken to be $\theta=7/64$.
\end{thm}
\noindent By considering such a mean value with a complex conjugate, we avoid some of the difficulties with Gauss sums and do not need to average over the modulus as in Liu's work. Since $L(1,f)$ is non-zero (see \cite[Lemma 5.9]{iwa}), we get
\begin{cor}
Fix $f$ a Hecke-Maass form for $\sltwo$. For every large enough prime $q$, there exists a primitive Dirichlet character $\chi$ of conductor $q$ such that the central $L$-values $L(\thalf, f \otimes \chi )$ and $L(\thalf, \chi)$ do not vanish. 
\end{cor}
\noindent Thus in our work, the modulus of the desired Dirichlet character is known. We have worked with prime moduli to minimize technical details. However given the power saving error term in (\ref{thmline}), we expect the same method would yield the non-vanishing result for any large integer $q$. 

Blomer, Fouvry, Kowalski, Michel, and Mili\'cevi\'c \cite{bloetal} have proven\footnote{At the time this paper was submitted for publication, \cite[Theorem 1.2]{bloetal} was conditional on the Ramanujan-Petersson conjecture, but now it is unconditional.}, conditionally on the Ramanujan-Petersson conjecture for $f$, the impressive asymptotic
\begin{align}
\label{blomean} \summ_{\chi \bmod q } L(\thalf, f \otimes \chi) \overline{L(\thalf, \chi)}^2 =  \frac{(q-2)L(1,f)^2}{\zeta(2)} +O(q^{-\delta}),
\end{align}
for $q$ prime and some $\delta>0$. Although the mean value we consider is more modest, it is not clear whether (\ref{thmline}) can be expected from (\ref{blomean}). Firstly, our asymptotic is unconditional and secondly, on a more technical level, the approximate functional equations needed for (\ref{blomean}) are set up without Gauss sums \cite[section 3]{bloetal} while for (\ref{thmline}), Gauss sums cannot be avoided.

Finally we remark that the problem we are considering is analytic in nature. For comparision, we mention that Chinta \cite{chinta} has shown that for a fixed elliptic curve $E$ over $\mathbb{Q}$ and $q$ a large prime, the central value $L(1/2, E \otimes \chi)$ is non-zero for all but $O(q^{7/8+\epsilon})$ Dirichlet characters $\chi$ modulo $q$. Given this powerful result, the simultaneous non-vanishing problems discussed above are almost trivial when posed for weight two holomorphic Hecke-cusp forms. For Hecke-Maass forms however, the algebraic methods used in \cite{chinta} are not applicable.

\subsection{Acknowledgements} 
We are grateful to Sheng-Chi Liu for providing us with his preprint \cite{liu}. We thank the Department of Mathematics, Indian Institute of science, Bangalore and Texas A\&M University at Qatar, where this work was done, for their hospitalities. The first author acknowledges the financial support by UGC-DST India during this work. 

\section{Preliminaries}

We give the proof of Theorem \ref{mainthm} for even Hecke-Maass forms, the details for odd forms being entirely similar. Thus throughout, $q$ will denote a prime and $f$ an even Hecke-Maass cusp form for the full modular group with Laplacian eigenvalue $\tfrac{1}{4} + T_f^2$, where $T_f$ is real.  Let $\lambda_f(n)$ denote the eigenvalue of the $n$-th Hecke operator corresponding to $f$. For more details on Maass forms, we refer the reader to \cite{iwaniec}.

For an even Dirichlet character $\chi$ modulo $q$, we define the $L$-functions
\begin{align}
L(s,\chi) = \sum \limits_{n=1}^\infty \chi(n) n^{-s} ,\q L(s,f \otimes \chi) = \sum \limits_{n=1}^\infty \lambda_f(n) \chi(n) n^{-s}
\end{align}
for $\Re(s)>1$, with analytic continuation to entire functions having the functional equations
\begin{align}
\Lambda(s,\chi) = \frac{\tau(\overline{\chi})}{\sqrt{q}} \Lambda(1-s, \overline{\chi}), \q   \Lambda(s,f \otimes \chi) = \frac{ \tau(\chi)^2}{q}  \Lambda(1-s,f \otimes \overline{\chi}),
\end{align}
where
\begin{align}
\Lambda(s,\chi):= \left( \frac{q}{\pi} \right)^{s/2} \Gamma \left( \frac{s}{2} \right) L(s,\chi) , \q  \Lambda(s,f \otimes \chi):= \left( \frac{q}{\pi} \right) ^{s}  \Gamma \left( \frac{s+ i T_f}{2} \right)\Gamma \left( \frac{s- i T_f}{2} \right) L(s,f \otimes \chi).
\end{align}

We will use the convention that $\epsilon$ denotes an arbitrarily small positive constant, but not necessarily the same constant from one occurrence to the next. All implied constants may depend implicitly on $\epsilon$ and $f$.

\subsection{Orthogonality of characters} 

We will need the following basic identity. This and more on Dirichlet characters can be found in \cite{dav}.
\begin{lem}\label{ortho}
For $q$ prime and $(nm,q)=1$, we have
\begin{align}
\summ_{\chi \bmod q} \chi(n)\overline{\chi}(m) = \begin{cases} q-2 &\text{ if } n\equiv m \bmod q \\ -1 &\text{ otherwise.} \end{cases}
\end{align}
\end{lem}

\subsection{Approximate functional equations}  

We will need the following expressions for the central values of $L(s, \overline{\chi})$ and $L(s, f \otimes \chi)$. These can be derived in a standard way from \cite[Theorem 5.3]{iwa} and the functional equations of these $L$-functions.

\begin{lem} \label{approx1}
For $\chi$ an even primitive character of modulus $q$, we have
\begin{align}
\label{afe1} \overline{L(\thalf, \chi)} = L(\thalf, \overline{\chi}) = \sum_{m \ge 1} \frac{\overline{\chi}(m)}{\sqrt{m}} V_1\left( \frac{m}{\sqrt{q}} \right) + \frac{\tau(\overline{\chi})}{\sqrt{q}} \sum_{m \ge 1} \frac{\chi(m)}{\sqrt{m}} V_1\left( \frac{m}{\sqrt{q}} \right),
\end{align}
where for $x,c>0$,
\begin{align}
V_1(x) = \frac{1}{2\pi i} \int_{(c)} (\sqrt{\pi}x)^{-s}  \frac{\Gamma\left(\frac{2s+1}{4}\right)}{\Gamma\left(\frac{1}{4}\right)} \frac{ds}{s}
\end{align}
We have that the $\ell$-th derivative $V_1^{(\ell)}(x) \ll_{\ell,C} \min\{ 1, x^{-C-\ell} \}$ for any $C>0$, and $V_1(x)\sim 1$ for $x< q^{-\epsilon}$. Thus the sums in (\ref{afe1}) are essentially supported on $m<q^{1/2+\epsilon}$.
\end{lem}

\begin{lem} \label{approx2}
For $\chi$ an even primitive character of modulus $q$, we have 
\begin{align}
\label{afe2} L(\thalf, f\otimes \chi) = \sum_{n \ge 1} \frac{\lambda_f(n)\chi(n) }{\sqrt{n}} V_2 \left( \frac{n}{q} \right) + \frac{ \tau(\chi)^2}{q} \sum_{n \ge 1} \frac{\lambda_f(n) \overline{\chi}(n)}{\sqrt{n}} V_2\left(  \frac{n}{q} \right),
\end{align}
where for $x,c>0$,
\begin{align}
V_2(x) =\frac{1}{2\pi i} \int_{(c)} (\pi x)^{-s}   \frac{\Gamma\left(\frac{2s+1+i2T_f}{4}\right) \Gamma\left(\frac{2s+1-i2T_f}{4}\right)}{\Gamma\left(\frac{1+i2T_f}{4}\right) \Gamma\left(\frac{1-i2T_f}{4}\right)}  \frac{ds}{s}.
\end{align}
We have that the $\ell$-th derivative $V_2^{(\ell)}(x) \ll_{\ell,C} \min\{ 1, x^{-C-\ell} \}$ for any $C>0$, and $V_2(x)\sim 1$ for $x<q^{-\epsilon}$. Thus the sums in (\ref{afe2}) are essentially supported on $n<q^{1+\epsilon}$.
\end{lem}

\subsection{Sums of Fourier coefficients.}

The Ramanujan conjecture for the Fourier coefficients of $f$ is true on average, by Rankin-Selberg theory. Namely, we have
\begin{align}
\label{ramavg} \sum_{n<x} |\lambda_f(n)| \ll x^{1/2}\left( \sum_{n<x} |\lambda_f(n)|^2 \right)^{1/2} \ll x.
\end{align}
Individually, we have Kim and Sarnak's \cite{kimsar} bound 
\begin{align} 
\label{ram} |\lambda_f(n)|\ll n^{\theta+\epsilon},
\end{align}
where $\theta=7/64$.

We will also encounter sums of Fourier coefficients twisted by additive characters. In this context, let us recall the Voronoi summation formula, 
\begin{lem} \label{voronoi} \cite[Theorem 4.2]{godber}
Let $\psi$ be a fixed smooth function with compact support on the positive reals. Let $d,\overline{d} \in \mathbb{Z}$ with $(q,d) = 1$ and $d\overline{d}\equiv1 \pmod{q}$. Then 
\begin{align}
\displaystyle\sum_{n\ge 1}\lambda_f(n)e\left(\frac{n\overline{d}}{q}\right) \psi\left(\frac{n}{N}\right) =q \displaystyle\sum_{n\ge 1}\frac{\lambda_f(n)}{n}e\left(\frac{nd}{q}\right) \Psi_+\left(\frac{nN}{q^2}\right) + q \displaystyle\sum_{n\ge 1}\frac{\lambda_f(n)}{n}e\left(\frac{-nd}{q}\right) \Psi_-\left(\frac{nN}{q^2}\right),
\end{align}
where for $\sigma > -1$,
\begin{align}
\label{psidef} \Psi_\pm (x) = \frac{ 1}{2\pi i}\int_{(\sigma)}(\pi^2x)^{-s} G_\pm(s) \widetilde{\psi}(-s)ds,
\end{align}
$\widetilde{\psi}(s)$ is the Mellin transform of $\psi(x)$ and 
\begin{align}
2\pi G_\pm(s) = \frac{ \Gamma\left(\frac{1+s+iT_f}{2}\right) \Gamma\left(\frac{1+s-iT_f}{2}\right) }{ \Gamma\left(\frac{-s+iT_f}{2}\right) \Gamma\left(\frac{-s-iT_f}{2}\right)} \pm \frac{ \Gamma\left(\frac{1+s+iT_f +1}{2}\right) \Gamma\left(\frac{1+s-iT_f +1}{2}\right) }{ \Gamma\left(\frac{-s+iT_f +1}{2}\right) \Gamma\left(\frac{-s-iT_f +1}{2}\right)}.
\end{align}
\end{lem}
\noindent The following result says that Fourier coefficients are orthogonal to additive characters on average.
\begin{lem} \label{addtwist} \cite[Theorem 8.1]{iwaniec} For any real number $\alpha$, we have that
\begin{align}
\sum_{n\le N} \lambda_f(n)e(\alpha n)  \ll N^{1/2+\epsilon}.
\end{align}
The implied constant does not depend on $\alpha$.
\end{lem}

\section{Proof of Theorem \ref{mainthm}}

Using the approximate functional equations and by picking out even characters using the factor $\frac{\chi(-1)+1}{2}$ which equals 1 when $\chi(-1)=1$ and 0 otherwise, we have that the left hand side of (\ref{thmline}) equals
\begin{multline}
 \frac{1}{2}\summ_{\chi \bmod q} (\chi(-1)+1)  \Bigg(  \sum_{n \ge 1} \frac{\lambda_f(n)\chi(n) }{\sqrt{n}} V_2 \left( \frac{n}{q} \right) + \frac{i^k \tau(\chi)^2}{q} \sum_{n \ge 1} \frac{\lambda_f(n) \overline{\chi}(n)}{\sqrt{n}} V_2\left(  \frac{n}{q} \right) \Bigg)  \\
\times \left(  \sum_{m \ge 1} \frac{\overline{\chi}(m)}{\sqrt{m}} V_1\left( \frac{m}{\sqrt{q}} \right) + \frac{\tau(\overline{\chi})}{\sqrt{q}} \sum_{m \ge 1} \frac{\chi(m)}{\sqrt{m}} V_1\left( \frac{m}{\sqrt{q}} \right) \right) .
\end{multline}
which we write as 
\begin{align}
 \frac{1}{2}\summ_{\chi \bmod q} (\chi(-1)+1) (S_1+S_2)(S_3+S_4).
\end{align}
Multiplying out the summand above leads to several cross terms, which we now analyze one by one.

\subsection{Cross terms with no Gauss sum} In this section we consider
\begin{multline}
 \frac{1}{2}\summ_{\chi \bmod q} (\chi(-1)+1)S_1S_3 =  \frac{1}{2}\summ_{\chi \bmod q} \sum_{n,m \ge 1} \frac{\lambda_f(n)\chi(n) \overline{\chi}(m)}{\sqrt{nm}} V_1\left( \frac{m}{\sqrt{q}} \right) V_2 \left( \frac{n}{q} \right) \\ +  \frac{1}{2}\summ_{\chi \bmod q} \sum_{n,m \ge 1} \frac{\lambda_f(n)\chi(-n) \overline{\chi}(m)}{\sqrt{nm}} V_1\left( \frac{m}{\sqrt{q}} \right) V_2 \left( \frac{n}{q} \right).
\end{multline}
\begin{lem} \label{s1b}
\begin{align}
\summ_{\chi \bmod q} \sum_{n,m \ge 1} \frac{\lambda_f(n)\chi(-n) \overline{\chi}(m)}{\sqrt{nm}} V_1\left( \frac{m}{\sqrt{q}} \right) V_2 \left( \frac{n}{q} \right) \ll q^{3/4+\theta + \epsilon}.
\end{align}
\end{lem}
\proof 
By Lemma \ref{ortho} and (\ref{ram}), it suffices to bound by $q^{3/4+\theta + \epsilon}$ the sum
\begin{align}
\label{s1a} q \sum_{\substack{ n<q^{1+\epsilon}, m<q^{1/2+\epsilon}\\ n+m \equiv 0 \bmod q}} \frac{ |\lambda_f(n)|}{\sqrt{nm}} \ll q^{1+\theta + \epsilon} \sum_{\substack{ n<q^{1+\epsilon}, m<q^{1/2+\epsilon}\\ n+m \equiv 0 \bmod q}} \frac{ 1}{\sqrt{nm}}.
\end{align}
Writing $n = kq - m$, we see from the ranges of $n$ and $m$ that we must have $1\le k < q^{\epsilon}$. Thus (\ref{s1a}) is bounded by
\begin{align}
q^{71/64 + \epsilon}  \sum_{m<q^{1/2+\epsilon}} \frac{1}{\sqrt{qm}} \ll q^{3/4+\theta + \epsilon}.
\end{align}
\endproof
The next result gives the main term of Theorem \ref{mainthm}.
\begin{lem}
\begin{align}
\label{diag}  \frac{1}{2}\summ_{\chi \bmod q} \sum_{n,m \ge 1} \frac{\lambda_f(n)\chi(n) \overline{\chi}(m)}{\sqrt{nm}} V_1\left( \frac{m}{\sqrt{q}} \right) V_2 \left( \frac{n}{q} \right)=  \frac{q-2}{2} L(1,f) + O(q^{3/4+\theta+\epsilon})
\end{align}
\end{lem}
\proof By Lemma \ref{ortho} and (\ref{ramavg}), the left hand side of (\ref{diag}) equals
\begin{align}
\label{s1c} \frac{q-2}{2} \sum_{\substack{n,m \ge 1 \\ n\equiv m \bmod q \\ (nm,q)=1}} \frac{\lambda_f(n)}{\sqrt{nm}} V_1\left( \frac{m}{\sqrt{q}} \right) V_2 \left( \frac{n}{q} \right) + O(q^{3/4+\epsilon}).
\end{align}
By a similar argument as that used to prove Lemma \ref{s1b}, we see that up to an error of $O(q^{3/4+\theta + \epsilon})$, the main term in (\ref{s1c}) consists of those terms with $n=m$:
 \begin{align}
\frac{q-2}{2} \sum_{\substack{n \ge 1\\ (n,q)=1} } \frac{\lambda_f(n)}{n} V_1\left( \frac{n}{\sqrt{q}} \right) V_2 \left( \frac{n}{q} \right) = \frac{q-2}{2} \sum_{n \ge 1 } \frac{\lambda_f(n)}{n} V_1\left( \frac{n}{\sqrt{q}} \right) V_2 \left( \frac{n}{q} \right) + O(q^{-100}).
\end{align}
The equality above holds because $V_1 (\frac{n}{\sqrt{q}})$ is very small unless $n<q^{1/2+\epsilon}$, in which case $(n,q)=1$ is automatic. Using the definitions of $V_1$ and $V_2$, we get
\begin{multline}
\frac{q-2}{2} \sum_{n \ge 1 } \frac{\lambda_f(n)}{n} V_1\left( \frac{n}{\sqrt{q}} \right) V_2 \left( \frac{n}{q} \right) = \\\frac{q-2}{2} \frac{1}{(2\pi i )^2} \int_{(1)} \int_{(1)} \sum_{n \ge 1 }  \frac{\lambda_f(n)}{n^{1+s_1+s_2}} (q \pi)^{-s_1/2 - s_2} \frac{\Gamma\left(\frac{2s_1+1}{4}\right)}{\Gamma\left(\frac{1}{4}\right)} \frac{\Gamma\left(\frac{2s_2+1+i2T_f}{4}\right) \Gamma\left(\frac{2s_2+1-i2T_f}{4}\right)}{\Gamma\left(\frac{1+i2T_f}{4}\right) \Gamma\left(\frac{1-i2T_f}{4}\right)}  \frac{ds_1}{s_1} \frac{ds_2}{s_2}.
\end{multline}
The $n$-sum inside the integral equals $L(1+s_1+s_2, f)$. Shifting the lines of integration to $\Re(s_1)=\Re(s_2)=-1/2+\epsilon$, we pick up the residue $L(1,f)$ at $s_1=s_2=0$. The integral on the new lines can be bounded in a standard way by $q^{-1/4+\epsilon}$. This completes the proof.
\endproof

\subsection{Cross terms with one Gauss sum} In this section we consider
\begin{align}
\summ_{\chi \bmod q} (\chi(-1)+1)S_2 S_4 \ \ \ \text{  and  } \ \ \ \summ_{\chi \bmod q} (\chi(-1)+1)S_1 S_4.
\end{align}

\begin{lem}
\begin{align}
\label{rem1}&\summ_{\chi \bmod q} (\chi(-1)+1) \cdot \frac{ \tau(\chi)^2}{q} \sum_{n \ge 1} \frac{\lambda_f(n) \overline{\chi}(n)}{\sqrt{n}} V_2\left(  \frac{n}{q} \right) \cdot  \frac{\tau(\overline{\chi})}{\sqrt{q}} \sum_{m \ge 1} \frac{\chi(m)}{\sqrt{m}} V_1\left( \frac{m}{\sqrt{q}} \right) \ll q^{3/4+\epsilon},\\
\label{rem2}&\summ_{\chi \bmod q} (\chi(-1)+1) \cdot  \sum_{n \ge 1} \frac{\lambda_f(n) \chi (n)}{\sqrt{n}} V_2\left(  \frac{n}{q} \right) \cdot  \frac{\tau(\overline{\chi})}{\sqrt{q}} \sum_{m \ge 1} \frac{\chi(m)}{\sqrt{m}} V_1\left( \frac{m}{\sqrt{q}} \right) \ll q^{3/4+\epsilon}.
\end{align}
\end{lem}
\proof The proofs of (\ref{rem1}) and (\ref{rem2}) are similar, so we show only the former. For even primitive characters we have that $\tau(\chi)\tau(\overline{\chi})=q$. Thus to prove (\ref{rem1}) we need to bound by $ q^{3/4+\epsilon}$ the sum
\begin{align}
\label{rem1a} \frac{1}{q^{1/2}} \sum_{\substack{n,m \ge 1\\(nm,q)=1}}   \frac{\lambda_f(n)}{\sqrt{nm}}  V_1\left( \frac{m}{\sqrt{q}}\right)   V_2\left(  \frac{n}{q} \right)  \summ_{\chi \bmod q} \tau(\chi) \overline{\chi}(\pm n) \chi(m  ).
\end{align}
Writing $\tau(\chi) =  \summ_{a \bmod q} \chi(a) e(a/q)$ , the innermost sum of (\ref{rem1a}) equals
\begin{align}
\summ_{\chi \bmod q} \tau(\chi) \overline{\chi}(\pm n) \chi(m  ) =  \summ_{a \bmod q} e\left(\frac{a}{q}\right) \summ_{\chi \bmod q} \overline{\chi}(\pm n) \chi(m a ).
\end{align}
Now using Lemma \ref{ortho}, we have that (\ref{rem1a}) equals
\begin{align}
\label{rem1b} &\frac{q-1}{q^{1/2}} \sum_{\substack{m \ge 1\\(m,q)=1}}  \frac{1}{\sqrt{m}} V_1\left( \frac{m}{\sqrt{q}}\right) \sum_{\substack{n\ge 1\\(n,q)=1}}  \frac{ \lambda_f(n) }{\sqrt{n}} e\left( \frac{\pm n \overline{m}}{q} \right) V_2\left(  \frac{n}{q} \right)  \\
\label{rem1c} &-\frac{1}{q^{1/2}} \sum_{\substack{m \ge 1\\(m,q)=1}}  \frac{1}{\sqrt{m}}  V_1\left( \frac{m}{\sqrt{q}}\right)  \sum_{\substack{n\ge 1\\(n,q)=1}}  \frac{\lambda_f(n)}{\sqrt{n}}  V_2\left(  \frac{n}{q} \right) \summ_{a \bmod q} e\left(\frac{a}{q}\right).
\end{align}
The innermost $a$-sum of (\ref{rem1c}) is a Ramanujan sum which equals $-1$, so that (\ref{rem1c}) is trivially bounded by $q^{1/4+\epsilon}$.  As for (\ref{rem1b}), removing the condition $(n,q)=1$ and using (\ref{ram}), we have that it equals
\begin{align}
\label{rem1d} \frac{q-1}{q^{1/2}} \sum_{\substack{m \ge 1\\(m,q)=1}}  \frac{1}{\sqrt{m}}  V_1\left( \frac{m}{\sqrt{q}}\right)   \sum_{n\ge 1}   \frac{\lambda_f(n)}{\sqrt{n}}  e\left( \frac{\pm n \overline{m}}{q} \right) V_2\left(  \frac{n}{q} \right) + O(q^{1/4+\theta+\epsilon}).
\end{align}
Now by Lemma \ref{addtwist} and partial summation, we get that the $n$-sum in \eqref{rem1d} is bounded by $q^\epsilon$, the $m$-sum by $q^{1/4}$. Hence \eqref{rem1d}   is bounded by $q^{3/4+\epsilon}$.
\endproof

\subsection{Cross terms with two Gauss sums}
In this section we consider
\begin{align}
\summ_{\chi \bmod q} (\chi(-1)+1)S_1 S_3.
\end{align}

\begin{lem}
\begin{align}
\summ_{\chi \bmod q} (\chi(-1)+1)  \cdot \frac{\tau(\chi)^2}{q} \sum_{n\ge 1}  \frac{\lambda_f(n) \overline{\chi} (n)}{\sqrt{nm}} V_2\left( \frac{n}{q}\right) \cdot  \sum_{m\ge 1} \frac{ \overline{\chi} ( m)}{\sqrt{m}}  V_1\left(\frac{m}{\sqrt{q}}\right) \ll q^{7/8+\theta+\epsilon}
\end{align}
\end{lem}
\proof
By taking a smooth partition of unity, we may consider the $n$ and $m$ sums in dyadic intervals. Thus it suffices to bound by $q^{7/8+\theta+\epsilon}$ the sum
\begin{align}
\label{2gauss} \frac{1}{(MN)^{1/2}} \summ_{\chi \bmod q} \frac{\tau(\chi)^2}{q} \sum_{n,m\ge 1} \lambda_f(n) \overline{\chi} (\pm nm) V_1\left(\frac{m}{\sqrt{q}}\right) W_1\left( \frac{m}{M} \right) V_2\left( \frac{n}{q}\right) W_2\left( \frac{n}{N} \right),
\end{align}
for any fixed smooth functions $W_i(x)$ supported on $x\in [1,2]$, $N<q^{1+\epsilon}$ and $M<q^{1/2+\epsilon}$. Writing $\tau(\chi) =  \summ_{a \bmod q} \chi(a) e(a/q)$, we have that (\ref{2gauss}) equals
\begin{multline}
\label{off} 
\frac{1}{q(MN)^{1/2}} \sum_{\substack{n,m\ge 1\\ (nm,q)=1}}  \lambda_f(n) V_1\left(\frac{m}{\sqrt{q}}\right) W_1\left( \frac{m}{M} \right) V_2\left( \frac{n}{q}\right) W_2\left( \frac{n}{N} \right) \\
\cdot \summ_{a,b \bmod q} e\left( \frac{a+b}{q}\right) \summ_{\chi \bmod q}  \overline{\chi}(\pm nm)\chi(ab).
\end{multline}
The innermost sum can be evaluated using Lemma \ref{ortho}. Thus (\ref{off}) equals
\begin{align}
\label{off2a} 
&\frac{q-1}{q(MN)^{1/2}} \sum_{\substack{n,m\ge 1\\ (nm,q)=1}}  \lambda_f(n) V_1\left(\frac{m}{\sqrt{q}}\right) W_1\left( \frac{m}{M} \right) V_2\left( \frac{n}{q}\right) W_2\left( \frac{n}{N} \right) 
\summ_{a \bmod q} e\left( \frac{a\pm nm\overline{a}}{q}\right)\\
\label{off2b} 
&-\frac{1}{q(MN)^{1/2}} \sum_{\substack{n,m\ge 1\\ (nm,q)=1}}  \lambda_f(n) V_1\left(\frac{m}{\sqrt{q}}\right) W_1\left( \frac{m}{M} \right) V_2\left( \frac{n}{q}\right) W_2\left( \frac{n}{N} \right) 
\summ_{a,b \bmod q} e\left( \frac{a+b}{q}\right)
\end{align}
The innermost $a,b$-sum of (\ref{off2b}) is a product of two Ramanujan sums and equals 1. Thus (\ref{off2b}) is bounded absolutely by $q^{-1/4+\epsilon}.$ As for (\ref{off2a}), we proceed according to the sizes of $N$ and $M$.

{\bf Case I:} $N< q^{1/2}, M<q^{1/4}$.
We note that the innermost $a$-sum of (\ref{off2a}) is a Kloosterman sum, which is less than $2q^{1/2}$ by Weil's bound. Using this and (\ref{ramavg}), we get that (\ref{off2a}) is bounded by
\begin{align}
q^{\epsilon}(qMN)^{1/2}\ll q^{7/8+\epsilon}.
\end{align}

{\bf Case II:} $M\ge q^{1/4}$. We first note that the $m$-sum in (\ref{off2a}) equals
\begin{align}
 \sum_{\substack{m \ge 1\\ (m,q)=1}} & e\left( \frac{\pm nm\overline{a}}{q} \right)  V_1\left( \frac{m}{\sqrt{q}} \right) W_1\left(\frac{m}{M}\right)  =  \sum_{ m \ge 1} & e\left( \frac{\pm nm\overline{a}}{q} \right)  V_1\left( \frac{m}{\sqrt{q}} \right) W_1\left(\frac{m}{M}\right) +O(q^{-100}),
\end{align}
because $V_1( \frac{m}{\sqrt{q}} )$ is very small unless $m<q^{1/2+\epsilon}$, in which case $(m,q)=1$ is automatic.
By Poisson summation, we have that
\begin{align}
 \sum_{m \ge 1} & e\left( \frac{\pm nm\overline{a}}{q} \right)  V_1\left( \frac{m}{\sqrt{q}} \right) W_1\left(\frac{m}{M}\right) \\ 
 \nonumber &= \sum_{1\le b \le q} e\left( \frac{\pm nb\overline{a}}{q} \right) \sum_{m\ge 1}   V_1\left( \frac{b+mq}{\sqrt{q}} \right) W_1\left(\frac{b+mq}{M}\right)\\
\nonumber &= \sum_{1\le b \le q} e\left( \frac{\pm nb\overline{a}}{q} \right) \sum \limits_{m= - \infty}^\infty \intt V_1\left( \frac{b+xq}{\sqrt{q}} \right) W_1\left(\frac{b+xq}{M}\right) e(-mx) dx\\
\nonumber &= \frac{M}{q} \sum \limits_{m=-\infty}^\infty \sum_{b \bmod q}  e\left( \frac{b(\pm n \overline{a} -m)}{q} \right)   \intt V_1\left(\frac{yM}{\sqrt{q}}\right) W_1(y) e\left( \frac{-mMy}{q} \right) dy.
\end{align}
Repeated integration by parts shows that the integral above is bounded by $(q^{1+\epsilon}/|m|M)^{B}$ for any $B\ge 0$. Thus we may restrict the $m$-sum in the last line to $|m|<q^{1+\epsilon}/M$, up to an error of $q^{-100}$ say. The $b$-sum equals $q$ if $\pm n \overline{a} \equiv m \bmod q$, and $0$ otherwise. Thus (\ref{off2a}) is bounded by
\begin{align}
& \frac{M^{1/2}}{N^{1/2}} \sum_{\substack{0< |m| < q^{1+\epsilon}/M \\ (m,q)=1}}\left| \sum_{\substack{n\ge 1\\ (n,q)=1}}  \lambda_f(n) e\left(\frac{\pm n\overline{m}}{q}\right) V_2\left( \frac{n}{q}\right) W_2\left( \frac{n}{N} \right)\right| + O(q^{-100}) \label{inter} \\
 &=\frac{M^{1/2}}{N^{1/2}} \sum_{\substack{0< |m| < q^{1+\epsilon}/M\\ (m,q)=1}}\left| \sum_{n\ge 1}  \lambda_f(n) e\left(\frac{\pm n\overline{m}}{q}\right) V_2\left( \frac{n}{q}\right) W_2\left( \frac{n}{N} \right)\right| +O(q^{1/4+\theta+\epsilon}). \label{case2}
\end{align}
Above, \eqref{case2} was obtained from \eqref{inter} by bounding absolutely the contributing of the integers $n$ divisible by $q$, of which there are at most $q^{\epsilon}$. By Lemma \ref{addtwist} and partial summation, we find that the $n$-sum in \eqref{case2} is bounded by $N^{1/2 + \epsilon}$ and so \eqref{case2} is bounded by
\begin{align}
\frac{q^{1+\epsilon}}{M^{1/2}} \ll q^{7/8+\epsilon}.
\end{align}

{\bf Case III:} $N\ge q^{1/2}, M<q^{1/4}$. By removing the condition $(n,q)=1$, we note that (\ref{off2a}) equals
\begin{multline}
\label{off2a1} \frac{q-1}{q(MN)^{1/2}} \sum_{\substack{m\ge 1\\ (m,q)=1}} V_1\left(\frac{m}{\sqrt{q}}\right) W_1\left( \frac{m}{M} \right) \summ_{a \bmod q} e\left( \frac{a}{q}\right)   \sum_{n\ge 1}  \lambda_f(n) V_2\left( \frac{n}{q}\right) W_2\left( \frac{n}{N} \right) e\left( \frac{\pm nm\overline{a}}{q}\right) \\+O(q^{-1/4+\theta+\epsilon}).
\end{multline}
Applying Voronoi summation to the innermost $n$-sum, we get that (\ref{off2a1}) equals
\begin{multline}
\label{nsum} \frac{q-1}{q(MN)^{1/2}} \sum_{\substack{m\ge 1\\ (m,q)=1}}  V_1\left(\frac{m}{\sqrt{q}}\right) W_1\left( \frac{m}{M} \right)\summ_{a \bmod q} e\left( \frac{a}{q}\right)  
\cdot q\sum_{n\ge 1} \frac{\lambda_f(n)}{n} e\left(  \frac{\mp n\overline{m}a}{q} \right) \Psi_- \left( \frac{nN}{q^2} \right) \\+ O(q^{-1/4+\theta+\epsilon})
\end{multline}
plus a similar sum involving $\Psi_+$, where 
\begin{align}
\label{psi} \Psi_\pm(x) =  \frac{1}{2\pi i}  \int_{(0)}(\pi^2x)^{-s} G_\pm(s) \int_{0}^{\infty} V_2(tN/q) W_2(t) t^{-s-1} dt \ ds.
\end{align}
The $a$-sum in (\ref{nsum}) equals $q-1$ if $n\equiv \mp m \bmod q$ and $-1$ otherwise, so that (\ref{nsum}) equals
\begin{align}
\label{vor1} & \frac{q(q-1)}{(MN)^{1/2}} \sum_{\substack{m\ge 1\\ (m,q)=1}}  V_1\left(\frac{m}{\sqrt{q}}\right) W_1\left( \frac{m}{M} \right)
\sum_{\substack{n\ge 1\\ n\equiv \mp m \bmod q}} \frac{\lambda_f(n)}{n} \Psi_- \left( \frac{nN}{q^2} \right)  \\
\label{vor2} &- \frac{(q-1)}{(MN)^{1/2}} \sum_{\substack{m\ge 1\\ (m,q)=1}}  V_1\left(\frac{m}{\sqrt{q}}\right) W_1\left( \frac{m}{M} \right)
\sum_{n\ge 1} \frac{\lambda_f(n)}{n} \Psi_- \left( \frac{nN}{q^2} \right) + O(q^{-1/4+\theta+\epsilon}).
\end{align}
We now explain how to truncate the $n$-sum. We first note that
\begin{align}
G_\pm(s) \int_{0}^{\infty} V_2(tN/q) W_2(t) t^{-s-1} dt  \ll (1+|s|)^{2\Re(s)+1} \left(\frac{q^\epsilon}{1+|s|}\right)^B,
\end{align}
by Stirling's approximation for the gamma function and by integrating by parts several times the $t$-integral, for any $B\ge 0$, where the implied constant depends on $\Re(s)$, $B$ and of course $f$. Using this estimate, by shifting the line of integration in (\ref{psi}) right to $\Re(s) = C$, we have that $\Psi_\pm(x)\ll q^{\epsilon}x^{-C}$ for any $C>0$. Thus we may restrict (\ref{vor1}) and (\ref{vor2}) to $n<q^{2+\epsilon}/N$, up to an error of $O(q^{-100})$ say. Also, we may shift the line of integration in (\ref{psi}) left to $\Re(s)=-1+\epsilon$ to get that that $\Psi_\pm(x)\ll q^{\epsilon} x$. 

Now restricting to $n<q^{2+\epsilon}/N$ and bounding absolutely we find that that (\ref{vor1}) is bounded by
\begin{align}
\frac{q^{2+\epsilon}}{(MN)^{1/2}} \sum_{m<  M} \sum_{\substack{n< q^{2+\epsilon}/N \\ n\equiv \mp m \bmod q}} \frac{|\lambda_f(n)|}{n} \frac{nN}{q^2} \ll \frac{q^{1+\theta+\epsilon}M^{1/2}}{N^{1/2}} \ll q^{7/8+\theta+\epsilon}.
\end{align}
The same bound holds for (\ref{vor2}). Since the sum invlolving $\Psi_+$ can be treated in exactly the same way, this completes the proof.

\bibliographystyle{amsplain}

\bibliography{simult-nonv}

\end{document}